\documentclass[11pt]{article}

\title{\sc A short proof that ${\mathcal B}(L_1)$ is not amenable}
\author{\sc Yemon Choi}

\date{12th October 2020}

\usepackage[left=30mm,right=30mm,top=30mm,bottom=30mm]{geometry}

\usepackage{sectsty}
\allsectionsfont{\sffamily}
\usepackage[colorlinks=true,citecolor=blue,linkcolor=black]{hyperref}

\input{BL1namac}

\begin{document}
\maketitle

\begin{abstract}
Non-amenability of ${\mathcal B}(E)$ has been surprisingly difficult to prove for the classical Banach spaces, but is now known for $E= \ell_p$ and $E=L_p$ for all $1\leq p<\infty$. However, the arguments are rather indirect: the proof for $L_1$ goes via non-amenability of $\ell^\infty({\mathcal K}(\ell_1))$ and a transference principle developed by Daws and Runde (Studia Math., 2010).

In this note, we provide a short proof that ${\mathcal B}(L_1)$ and some of its subalgebras are non-amenable, which completely bypasses all of this machinery. Our approach is based on classical properties of the ideal of representable operators on $L_1$, and shows that ${\mathcal B}(L_1)$ is not even approximately amenable.

\bigskip
\noindent
Keywords: amenable Banach algebras, Banach spaces, operator ideals, representable operators.

\bigskip
\noindent
MSC 2020: 46H10, 47L10 (primary); 46B22, 46G10 (secondary)

\end{abstract}

\begin{section}{Introduction}
Throughout this paper: all algebras are associative and taken over the field $\Cplx$, but they need not have identity elements.

The Wedderburn structure theorem implies, with hindsight, that a finite-dimensional algebra with homological dimension zero is isomorphic to a sum of full matrix algebras. Amenability for Banach algebras, introduced in B. E. Johnson's seminal work \cite{BEJ_CIBA}, may be thought of as a weakened version of having homological dimension zero,
 and the two notions coincide for finite-dimensional Banach algebras.
In particular, finite sums of full matrix algebras are amenable, while the algebra of $2\times 2$ upper triangular matrices is not.

It is therefore natural to ask for which infinite-dimensional Banach spaces $E$ the algebra $\Bdd(E)$ is amenable. It was soon recognized that for most $E$ the answer should be negative, but that proving this for specific natural $E$ could be very hard.
While the Hilbertian case was known to follow very indirectly from deep results on ${\rm C}^*$-algebras, no progress was made on the other classical Banach spaces until C. J. Read's breakthrough result that $\Bdd(\ell_1)$ is non-amenable~\cite{Read_Bl1-nonamen}.
His proof was simplified by G. Pisier \cite{Pisier-on-Read}, and
 N. Ozawa subsequently provided a unified proof of non-amenability of $\Bdd(\ell_1)$, $\Bdd(\ell_2)$ and $\Bdd(c_0)$~\cite{Oza_Blp-nonamen}. Further historical details can be found in
V.~Runde's survey article \cite{Runde_BE-survey}, or in the introduction of his companion paper~\cite{Run_JAMS}.

The paper \cite{Run_JAMS} contains the strongest
 general results thus far on non-amenability of $\Bdd(E)$; among other things, it establishes the non-amenability of $\Bdd(\ell_p)$ and $B(L_p)$ for all $p\in (1,\infty)$. A key ingredient in the proof is the following ``transference principle'' developed by M. Daws and V. Runde in \cite{DawsRunde}:
\begin{itemize}
\item[--] if $F$ is a Banach space, amenability of $\Bdd(\ell_p(F))$ implies amenability of $\ell^\infty(\cK(\ell_p(F)))$;
\item[--] if $E$ is an infinite-dimensional $\cL_p$-space in the sense of Lindenstrauss and Pe{\l}czy\'nski, then amenability of $\ell^\infty(\cK(E))$ implies amenability of $\ell^\infty(\cK(X))$ for \emph{every} $\cL_p$-space $X$.
\end{itemize}
Although the case of $L_1$ was not resolved in \cite{Run_JAMS}, the transference principle remains valid for $p=1$ (see \cite[Theorems 1.2 and 4.3]{DawsRunde}), and so amenability of $\Bdd(L_1)\iso\Bdd(\ell_1(L_1))$ would imply amenability of $\ell^\infty(\cK(\ell_1))$. Therefore it suffices to prove that the latter algebra is non-amenable, and this was recently demonstrated in the PhD thesis of E.~Aldabbas~\cite{aldabbas}; as in \cite{Run_JAMS}, crucial use is made of a technical innovation from~\cite{Oza_Blp-nonamen}, which was itself inspired by the arguments of \cite{Read_Bl1-nonamen}. At the time of writing, the proof from \cite{aldabbas} has not been published.

Thus, although non-amenability of $\Bdd(L_1)$ is now known, the existing proof is both indirect (going via $\ell^\infty(\cK(\ell_1))$) and technically complicated (relying on ``Ozawa's lemma'' as formulated in \cite{Run_JAMS}). The purpose of this note is to show that non-amenability of $\Bdd(L_1)$ can be proved very quickly by studying a particular closed ideal $\fR\lhd \Bdd(L_1)$, without any need for the transference principle or the ideas in~\cite{Oza_Blp-nonamen}. Our method actually proves slightly more: if $A\subseteq\Bdd(L_1)$ is a closed subalgebra that contains~$\fR$, then $A$ is not even approximately amenable in the sense of \cite{GhL_genam1,GhLZ_genam2}.

The ideal $\fR$ occurs very naturally in the study of operators on $L_1$, and is related to a factorization result of D. R. Lewis and C. Stegall. Thus our approach has a rather different flavour from the combinatorial arguments in \cite{Read_Bl1-nonamen} and \cite{Oza_Blp-nonamen}, and we hope that this alternative perspective could be useful for studying the non-amenability problem for $\Bdd(E)$ on other Banach lattices.
\end{section}

\begin{section}{Background and preliminaries}

\begin{subsection}{Notation and other conventions}
Most of our conventions for notation and terminology are either standard in the literature or clear from context. However, to make our paper more accessible, we have endeavoured to provide precise references for various ``well-known'' or ``standard'' facts about Banach spaces.

The term ``operator'' is synonymous with ``bounded linear map'' although we shall sometimes refer to ``bounded operators'' just for emphasis.
For a Banach space $E$, $\Bdd(E)$ denotes the algebra of bounded operators on $E$, while $\cK(E)$ denotes the algebra of compact operators on~$E$.

The projective tensor product of Banach spaces $E$ and $F$ is denoted by $E\ptp F$. For our purposes, it can be characterized as by the following universal property: whenever $E$, $F$ and $G$ are Banach spaces and $\beta:E\times F\to G$ is a bounded bilinear map, there is a unique bounded linear map $f:E\ptp F \to G$ that satisfies $f(x\tp y) = \beta(x,y)$ for all $x\in E$ and $y\in F$. Moreover, $\norm{f}=\norm{\beta}$.
Note also that for any $x\in E$ and $y\in F$ we have $\norm{x\tp y}_{E\ptp F}=\norm{x}\norm{y}$.

Given a measure space $(\Om,\Sigma,\mu)$ we abbreviate $L_p(\Om,\Sigma,\mu)$ to $L_p(\Om)$. In the case of $[0,1]$ with the Borel $\sigma$-algebra and Lebesgue measure, we simply write~$L_p$; in the case of $\Nat$ with the discrete $\sigma$-algebra and counting measure, we simply write~$\ell_p$.
If $A$ is a Borel subset of $[0,1]$ then $|A|$ denotes its Lebesgue measure.

The background we need concerning Banach-space valued integration can be found in any source that defines the Bochner integral over a finite measure space. By a slight abuse of terminology, we say that a function from $[0,1]$ to a Banach space $X$ is \dt{strongly measurable} if it is strongly measurable with respect to the Borel $\sigma$-algebra of $[0,1]$.
\end{subsection}

\begin{subsection}{Amenable Banach algebras}
Since this paper is intended for a general rather than a specialist audience, we use this section to record some basic definitions and examples from the literature on amenability of Banach algebras, in order to supply some context for the main result.

The following definition of amenability is not the original one given by B. E. Johnson in \cite{BEJ_CIBA}, but is a standard equivalent formulation that is often more useful or more suggestive.

\begin{dfn}
Let $A$ be a Banach algebra and define $\pi:A\ptp A\to A$ by $\pi(a\tp b)=ab$. An \dt{approximate diagonal for~$A$} is a net $(d_\alpha)$ in $A\ptp A$ such that, for each $a\in A$, we have
\[
\norm{a\cdot d_\alpha - d_\alpha \cdot a}_{A\ptp A} \to 0
\quad\text{and}\quad
\norm{a\pi(d_\alpha)-a}_A \to 0
\quad\text{as $\alpha\to\infty$.}
\]
If the net $(d_\alpha)$ is bounded then we call it a \dt{bounded approximate diagonal}. A Banach algebra possessing a bounded approximate diagonal is said to be \dt{amenable}.
\end{dfn}

Note that if $A$ is finite-dimensional and amenable, taking a cluster point of the net $(d_\alpha)$ yields an element $\Delta\in A\tp A$ such that $\pi(\Delta)=1_A$ and $a\cdot\Delta=\Delta\cdot a$. In (non-Banach) homological algebra such a $\Delta$ is known as a \dt{separating idempotent} for $A$ and serves as an explicit witness that $A$ has homological dimension~zero.

It is well known that finite-dimensional matrix algebras $M_n(\Cplx)\equiv \Bdd(\Cplx^n)$ have homological dimension-zero, and that an explicit separating idempotent for $M_n(\Cplx)$ is
\[
\Delta = \frac{1}{n}\sum_{i,j=1}^n E_{ij} \tp E_{ji} \;.
\]
By an averaging argument, $\Delta$ can be written as an absolutely convex combination of tensors of the form $x\tp x^{-1}$ where $x$ is a signed permutation matrix.
For details see the proof of \cite[Prop.~3.2]{GJW_amen-cpct}.
It follows that if $A=\Bdd(\ell_p^n)$ for any $1\leq p \leq\infty$, $\norm{\Delta}_{A\ptp A}=1$.
From this, a routine exhaustion argument allows one to construct an explicit bounded approximate diagonal for $\cK(\ell_p)$ when $1\leq p <\infty$ and $\cK(c_0)$.

Using a more abstract version of this idea, the paper \cite{GJW_amen-cpct} developed a condition on a given Banach space $E$, called Property $(\bbA)$, which is sufficient for amenability of~$\cK(E)$. Property $(\bbA)$ is studied in detail in that paper; while it is not known to hold for all $\cL_p$-spaces, it does hold for all $L_p(\mu)$-spaces ($1\leq p\leq\infty$) and their preduals (see \cite[Theorem 4.7]{GJW_amen-cpct} and \cite[Theorem 4.3]{GJW_amen-cpct}).

As mentioned in the introduction, the expectation is that for most $E$ the Banach algebra $\Bdd(E)$ is in some sense too large to be amenable. However, S. A. Argyros and R. Haydon constructed in \cite{ArgHay_S+C} an infinite-dimensional space $X$ such that every bounded operator on $X$ is a compact perturbation of a multiple of the identity, solving one of the major foundational problems of the subject. The nature of their construction also ensures that $X^*$ is isomorphic to $\ell_1$; thus $X$ has Property $(\bbA)$, and so $\cK(X)$ is amenable. Since unitizations of amenable Banach algebras are amenable, it follows that $\Bdd(X)=\Cplx I + \cK(X)$ is amenable (\cite[Prop.~10.6]{ArgHay_S+C}).

\end{subsection}

\begin{subsection}{Preliminary results needed for our paper}

The second condition in the definition of a (bounded) approximate diagonal says that the net $(\pi(d_\alpha))$ is a \dt{(bounded) right approximate identity} for~$A$. There is a corresponding notion of a (bounded) left approximate identity.
Crucially, amenability of a Banach algebra not only ensures bounded left and right approximate identities in the algebra itself, but also in some of its closed ideals. The following lemma follows from standard results in the theory of amenable Banach algebras and their bimodules.
For instance, it is an immediate corollary of \cite[Theorem~3.7]{CuLo_amen}.
\end{subsection}

\begin{lem}\label{l:both sides}
Let $A$ be an amenable Banach algebra and let $J$ be a closed, $2$-sided ideal in $A$ which has a bounded right approximate identity. Then $J$ has a bounded left approximate identity.
\end{lem}

From this we immediately deduce the following result.

\begin{cor}[An obstruction to amenability]\label{c:amen-obstruction}
Let $A$ be a Banach algebra and $J$ a closed $2$-sided ideal.
Suppose that
\begin{itemize}
\item $J$ has a bounded right approximate identity;
\item there exists $x_0\in J$ such that $x_0 \notin \overline{Jx_0}$.
\end{itemize}
Then $A$ is not amenable.
\end{cor}


In the next section, we will introduce a particular closed ideal in $\Bdd(L_1)$ and show that it satisfies both conditions in Corollary \ref{c:amen-obstruction}.

\end{section}

\begin{section}{The key ideal in $\Bdd(L_1)$}
\label{s:representable}
Given a Banach space $E$, an operator $T:L_1\to E$ is said to be \dt{representable} if there exists a bounded, strongly measurable function $h_T:[0,1]\to E$ such that
\[ T(f) = \int_0^1 f(s)h_T(s)\,ds \qquad\text{for all $f\in L_1$}
\]
where the right-hand side is interpreted as an $E$-valued Bochner integral.
(In some sources, such as \cite{Ros75_survey}, the terminology ``differentiable'' is used instead of ``representable''.)

Note that if such an $h_T$ exists, we have $\norm{T} \leq \norm{h_T}_\infty$ by basic properties of the Bochner integral; one can show that equality holds, although this is not needed for the present paper.

We denote by $\fR$ the set of all representable operators from $L_1$ to itself.
It follows easily from the definitions that $\fR$ is a closed left ideal in $\Bdd(L_1)$. Therefore, to show that it is also a right ideal, it suffices to prove that $TS\in \fR$ for all $T\in \fR$ and all $S\in \Bdd(L_1)$. This is not obvious from the definition, but is an immediate consequence of the next result  which is due to D. R. Lewis and C. Stegall.

\begin{thm}[Lewis--Stegall]\label{t:lewis-stegall}
Let $E$ be a Banach space and let $T\in\Bdd(L_1,E)$. Then $T$ is representable if and only if it factors (boundedly) through $\ell_1$.
\end{thm}

For a direct and relatively self-contained proof, which only needs the basic properties of strongly measurable $E$-valued functions, 
see~\cite[Appendix~C, \S4]{DefFlor}.
(Alternative sources  are \cite[Theorem~A3]{Ros75_survey}, \cite[Chapter III, Theorem 1.8]{DU_VMbook} or \cite[Prop.\ 5.36]{Ryan_TP}.)

\begin{rem}
The original paper \cite{LewSteg73} does not make use of the perspective of Bochner integrals and vector-valued $L_p$-spaces.
Indeed, while Theorem~\ref{t:lewis-stegall} was known at the time to follow from the techniques in \cite{LewSteg73}, the result itself is never explicitly stated there, although some version of it appears {\it en passant} in the proof of \cite[Theorem~1]{LewSteg73}. C.f.\ the remarks in \cite[Appendix A]{Ros75_survey}.
\end{rem}


The first part of the next result is well-known to Banach space theorists, although we are not aware of a reference.

\begin{prop}\label{p:no LAI}
There exists $T_0\in\fR\setminus\cK(L_1)$ such that $ST_0\in \cK(L_1)$ for all $S\in \fR$. In particular, $\fR$ does not have any left approximate identity (bounded or otherwise).
\end{prop}

\begin{proof}
Let $S\in\fR$. By the Lewis--Stegall theorem, $S$ factors through $\ell_1$, and hence it maps weakly convergent sequences to norm convergent sequences (since $\ell_1$ has the Schur property). By the Eberlein--\v{S}mulian theorem, it follows that $S$ maps relatively weakly compact subsets of $L_1$ to totally bounded subsets of $L_1$. Moreover, every weakly compact operator on $L_1$ is representable (see~e.g.\ \cite[Section~C5]{DefFlor}, \cite[Chapter III, Lemma 2.9]{DU_VMbook} or \cite[Prop.\ 5.40]{Ryan_TP}). It therefore suffices to choose any $T_0\in\Bdd(L_1)$ which is weakly compact but not compact.
\end{proof}

There are various indirect ways to show the existence of weakly compact non-compact operators on $L_1$. We describe one easy and explicit construction for the reader's convenience.

\begin{eg}\label{eg:explicit}
Fix a partition of $(0,1]$ into countably many subsets with strictly positive measure (e.g.\ intervals $(2^{-n},2^{1-n}]$ for $n\in\Nat$) and let $P:L_1\to \ell_1$ be the associated conditional expectation.
Let $\iota_{1,2} : \ell_1 \to \ell_2$ be the canonical embedding;
and fix an injection $j:\ell_2\to L_1$ with closed range (for instance, using Rademacher functions).
Then $T_0\defeq j\iota_{1,2}P$ is non-compact, since $P$ is an open mapping, $\iota_{1,2}$ is non-compact, and $j$ is bounded below.
On the other hand, $T_0$ is weakly compact since it factors through $\ell_2$. 
Note also that by Pitt's theorem we get a direct proof that $ST_0\in\cK(L_1)$ for all $S:L_1\to\ell_1$, without requiring the Lewis--Stegall theorem or the fact that weakly compact operators on $L_1$ are representable.
\end{eg}

W. B. Johnson has informed the author that in forthcoming work with N. C. Phillips and G. Schechtman, they show that for $1\leq p <\infty$ the only closed ideal in $\Bdd(L_p)$ with a bounded left approximate identity is $\cK(L_p)$. In the same work, they also establish the following result, which is the key ingredient needed for the present paper.

\begin{prop}[Johnson--Phillips--Schechtman]\label{p:WBJ_BRAI}
$\fR$ has a bounded right approximate identity. Moreover, we can choose this net to consist of norm-one idempotents.
\end{prop}

Since the work of Johnson--Phillips--Schechtman is still unpublished at time of writing, we include a self-contained proof of Proposition~\ref{p:WBJ_BRAI}.
The argument originally shown to the author by W. B. Johnson used ideas from \cite{LewSteg73} and some auxiliary results on $\cK(L_1)$. Our approach uses the perspective of vector-valued $L_\infty$-spaces, and is based on a suggestion of M.~Daws (personal communication).

\begin{lem}\label{l:countable}
Let $E$ be a Banach space and let $h:[0,1]\to E$ be strongly measurable.
For any $\veps>0$, there is a (strongly) measurable $h_\veps:[0,1]\to E$ whose range is \emph{countable} and which satisfies $\norm{h-h_\veps}_\infty \leq\veps$.
\end{lem}

We omit the proof of this lemma, which is a variation on the usual argument for scalar-valued functions. It is usually found in the literature as part of the proof of the Pettis measurability criterion (see e.g. \cite[Theorem B11]{DefFlor} or the proof of \cite[Prop.\ 2.15]{Ryan_TP}). For an explicit reference with a full proof, see \cite[Lemma 2.1.4]{HVVW_vol1}.

\begin{proof}[Proof of Proposition~\ref{p:WBJ_BRAI}]
Let $\fR_0$ be the set of operators $L_1\to L_1$ that are represented by bounded and \emph{countably-valued} measurable functions $[0,1]\to L_1$. Then $\fR_0$ is a left ideal in $\Bdd(L_1)$ and by Lemma~\ref{l:countable} it is dense in $\fR$. Hence, by a $3$-epsilon argument, it suffices to prove that $\fR_0$ has a bounded right approximate identity consisting of norm-one idempotents.

Given a partition of $[0,1]$ as a countable disjoint union of measurable subsets, $[0,1]=\bigsqcup_{n=1}^\infty A_n$, define a corresponding conditional expectation $P:L_1\to L_1$
by the formula
\begin{equation}\label{eq:cond-exp}
P(f)(t) = \frac{1}{|A_n|} \int_{A_n} f \qquad\text{if $t\in A_n$,}
\tag{$*$}
\end{equation}
with the convention that if $|A_n|=0$ we interpret $|A_n|^{-1}\int_{A_n} f$ as zero.
Then $P$ is a norm-one idempotent in $\Bdd(L_1)$, which belongs to $\fR_0$ since $P$ is represented by $h_P\defeq \sum_{n=1}^\infty |A_n|^{-1} 1_{A_n}$.

If $h:[0,1]\to L_1$ is constant on each $A_n$, with $h(A_n)=\{c_n\}$ say, then the operator $T\in\fR_0$ represented by $h$ satisfies
\[
TP(f)
= \sum_{n=1}^\infty \int_{A_n} h \cdot  (Pf)
= \sum_{n=1}^\infty c_n \int_{A_n} Pf
= \sum_{n=1}^\infty c_n \int_{A_n} f
= T(f)
\qquad\text(f\in L_1); \]
that is, $TP=T$.
Now, given $T_1,\dots, T_m \in\fR_0$, represented by bounded functions $h_1,\dots, h_m:[0,1]\to L_1$ respectively, note that there is a countable partition $[0,1]=\bigsqcup_{n=1}^\infty A_n$ such that each $h_j$ is constant on each $A_n$. Defining $P$ by the formula \eqref{eq:cond-exp}, the previous calculation now gives $T_jP=T_j$ for all $j=1,\dots, m$.

Therefore, if we order the set of countable partitions of $[0,1]$ by refinement, we obtain a net of norm-one idempotents in $\fR_0$, each having the form \eqref{eq:cond-exp}, which serves as a right approximate identity for $\fR_0$.
\end{proof}

Combining Proposition~\ref{p:no LAI} and Proposition~\ref{p:WBJ_BRAI}, we see that $\fR$ satisfies the conditions of Corollary~\ref{c:amen-obstruction}, and therefore $\Bdd(L_1)$ is not amenable. In fact, the corollary rules out amenability for every closed subalgebra $A\subseteq\Bdd(L_1)$ that contains $\fR$.

It is notable that for Proposition~\ref{p:no LAI}, the key feature of $\fR$ was that every $T\in\fR$ factors through $\ell_1$, while for Proposition~\ref{p:WBJ_BRAI} it seems vital to have the description in terms of representability by strongly measurable functions on~$[0,1]$.

\begin{rem}
In this section we chose to work with $\fR$ and its properties because it is an ideal with intrinsic interest, regardless of the application to non-amenability.
 One can bypass explicit mention of $\fR$ and extract a slightly more direct proof that $\Bdd(L_1)$ is non-amenable, by combining specific properties of the operator $T_0$ in Example~\ref{eg:explicit} with calculations in Appendix~\ref{app:both-sides-direct}.
However, this direct approach still seems to require the result that every operator $L_1\to\ell_1$ is representable (the ``easy direction'' of the Lewis--Stegall theorem), and so we do not include the details here.
\end{rem}

\end{section}

\begin{section}{Related examples and variations}
\label{s:loose ends}

Corollary~\ref{c:amen-obstruction} can be applied to prove non-amenability of $\Bdd(E)$ for some other Banach spaces~$E$, provided we make a left-right switch.
Since a Banach algebra $A$ is amenable if and only if the opposite algebra $A^{\rm op}$ is, Lemma~\ref{l:both sides} remains true when the words ``left'' or ``right'' are interchanged. Therefore, if a Banach algebra $A$ possesses a closed ideal $J$ that has a bounded left approximate identity but no bounded right appproximate idenity, $A$ cannot be amenable.

\begin{eg}
Let $E$ be a Banach space and let $\cA(E)$ denote the algebra of \dt{approximable operators} on~$E$; this is a closed ideal in $\Bdd(E)$. By results of N. Gr{\o}nb{\ae}k and G. A. Willis, $\cA(E)$ has a bounded left approximate identity if and only if $E$ has the bounded approximation property, but has a bounded right approximate identity if and only if $E^*$ has the bounded approximation property.
(See \cite[Theorem 2.4 and Theorem 3.3]{GroWil_CMB93}.)

Hence, by our previous remarks, if $E$ has the bounded approximation property and $E^*$ does not then $\Bdd(E)$ is non-amenable. This applies for instance when $E=\ell_2\ptp\ell_2$.
\end{eg}



It is natural to wonder if the techniques in this paper can be adapted to give an alternative proof of the non-amenability of $\Bdd(\ell_1)$. In this context, note that by \cite[Theorem 1.2]{DawsRunde}, amenability of $\Bdd(\ell_1)$ would imply amenability of $\ell^\infty(\Bdd(\ell_1))$ and hence amenability of any ultrapower ${\Bdd(\ell_1)}_{\mathcal U}$; such an ultrapower can be represented as an algebra of operators on some abstract $L$-space $E$, and if we can find an operator on $E$ analogous to the operator $T_0$ in Example~\ref{eg:explicit} then it may be possible to run similar arguments to the ones in this paper.
We leave this as a problem for possible future investigation.

We briefly comment on approximate amenability, although this was not the main focus of the present work. Given a Banach algebra $A$ let $A^\#$ denote its forced unitization. We say that $A$ is \dt{approximately amenable} if $A^\#$ has an approximate diagonal. This is not the original definition from \cite{GhL_genam1}; strictly speaking, what we have just defined is ``approximate contractibility'' of~$A$, but the two concepts were shown to coincide in \cite[Theorem 2.1]{GhLZ_genam2}.
By \cite[Corollary 2.4]{GhL_genam1}, one has an analogue of Lemma~\ref{l:both sides}:
\begin{quote}
if $A$ is approximately amenable and $J$ is a closed ideal in $A$ possessing a \emph{bounded} right approximate identity, then $J$ has a left approximate identity (not necessarily bounded).
\end{quote}
For an outline of a direct proof, see the appendix. From this result, we see that Corollary~\ref{c:amen-obstruction} remains valid if we weaken the hypothesis from amenability to approximate amenability. Hence, by the results of Section~\ref{s:representable},
every closed subalgebra of $\Bdd(L_1)$ which contains $\fR$ fails to be aproximately amenable.

\end{section}

\subsection*{Acknowledgements}
The question of whether the ideal $\fR$ has either a left or right bounded approximate identity emerged from discussions with Matt Daws and Jon Bannon about a completely different problem in von Neumann algebras, while the latter was visiting Lancaster University in 2019 as a Fulbright Visiting Scholar.

The author thanks Bence Horv\'ath and Niels Laustsen for spurring him to learn more about ideals in $\Bdd(L_p)$, Bill Johnson for answering his original questions about $\fR$ and clarifying some of the relevant history, Matt Daws for helpful discussions concerning Lebesgue--Bochner spaces, and Jon Bannon for reminding him that there must be better songs to sing than this.

\appendix

\begin{section}{A direct proof of Lemma~\ref{l:both sides}}
\label{app:both-sides-direct}
For sake of brevity, we justified the claim in Lemma~\ref{l:both sides} by appealing to more general results in \cite{CuLo_amen}. Specifically, we were invoking the following standard result:
\begin{quote}
if $A$ is an amenable Banach algebra and $J\lhd A$ is a closed ideal that is weakly complemented in $A$ as a Banach space, then $J$ has a bounded approximate identity.
\end{quote}
This result implies Lemma~\ref{l:both sides} because in any Banach algebra (regardless of amenability), a closed ideal with a bounded left-or-right approximate identity is weakly complemented.

The proof of the general result is somewhat abstract: one starts with a bounded linear projection from $A^*$ onto $J^\perp$, and then uses amenability to average this projection to an $A$-bimodule map, from which one extracts a left identity for $J^{**}$ equipped with the first Arens product.
It is therefore instructive to have a more direct proof of Lemma~\ref{l:both sides}, since this sheds more light on possible refinements of Corollary~\ref{c:amen-obstruction}.
We provide details below, since we have not seen such a proof written down explicitly. No novelty is claimed for the following arguments.

Let ${\rm FIN}(J)$ denote the set of finite subsets of $J$.
Our bounded left approximate identity will be indexed by 
${\rm FIN}(J)\times (0,\infty)$, given the following partial order: $(F,\veps)\preceq (F',\veps')$ if $F\subseteq F'$ and $\veps\geq \veps'$.
Thus, fix some $F\in {\rm FIN}(J)$ and $\veps>0$; it suffices to find $v\in J$ such that $\max_{x\in F} \norm{x-vx}<\veps$, and such that $\norm{v}$ is bounded above by a constant independent of $F$ and $\veps$.

The hypotheses of Lemma~\ref{l:both sides} ensure that for some constant $C>0$, $A$ has an approximate diagonal bounded in norm by $C$ and $J$ has a right approximate identity bounded in norm by~$C$. Let $\delta>0$ which will be chosen with hindsight to depend on $C$ and $\veps$.
Perturbing the bounded approximate diagonal slightly, we obtain $\Delta\in A\tp A$ with $\norm{\Delta}_{A\ptp A} \leq  C+1$ and
\begin{equation}\label{eq:approx-central}
\norm{x\cdot\Delta -\Delta\cdot x}_{A\ptp A} \leq\delta \quad\text{and}\quad \norm{x-x\pi(\Delta)}\leq \delta \quad\text{for all $x\in F$.}
\end{equation}
By definition of the projective tensor norm, we can assume that $\Delta=\sum_{i=1}^m a_i \tp b_i $ where $\sum_{i=1}^m \norm{a_i} \norm{b_i} \leq C+1$.
Since $\{xa_i \colon x\in F, 1\leq i\leq m\}$ is a finite subset of $J$, there exists some $f\in J$ with $\norm{f}\leq C$ and
\begin{equation}\label{eq:sandwich}
\norm{xa_i-xa_if}\leq \delta\norm{a_i} \quad\text{for all $x\in F$ and all $1\leq i\leq m$.}
\end{equation}

We put $v\defeq \sum_{i=1}^m a_i fb_i \in J$, which satisfies $\norm{v} \leq \norm{\Delta}_{A\ptp A}\norm{f} \leq C^2$. For each $x\in F$,
\begin{equation}\label{eq:three pieces}
\norm{ x- v x }
\leq
\norm{x - x\pi(\Delta)} + 
\norm{x \pi(\Delta)-x v} +
\norm{x v - vx} \;.
\end{equation}
The first term on the right-hand side is bounded above by $\delta$. The second term is bounded above by
\begin{equation}
\left\Vert x\sum_{i=1}^m a_ib_i - x\sum_{i=1}^m a_i fb_i\right\Vert
\leq \sum_{i=1}^m \norm{xa_ib_i -xa_i fb_i} \leq \sum_{i=1} \delta\norm{a_i}\norm{b_i} \leq \delta(C+1).
\end{equation}
To control the third term in \eqref{eq:three pieces}, note that the map $\theta: A\ptp A\to \Bdd(A)$ defined by $\theta(a\tp b)(z) =azb$ is contractive. Therefore, since $xv = \theta(x\cdot\Delta)(f)$ and $vx=\theta(\Delta\cdot x)(f)$,
\begin{equation}\label{eq:pincer trick}
\norm{xv-vx} =
\norm{ \theta(x\cdot\Delta - \Delta\cdot x)(f)}
\leq
{\norm{ x\cdot\Delta - \Delta\cdot x}}_{A\ptp A} \norm{f}
\leq \delta C.
\end{equation}
Hence $\norm{x-vx}\leq 2(C+1)\delta$; provided that we chose $\delta$ to ensure $2(C+1)\delta\leq\veps$, we have obtained the desired $v=v_{F,\veps}$.
This completes the proof of Lemma~\ref{l:both sides}.
\hfill$\Box$

\bigskip
We briefly indicate how one can adapt this argument to prove the ``approximately amenable version'' of Lemma~\ref{l:both sides} that was stated in Section~\ref{s:loose ends}. First, note that if $J$ is a closed ideal in $A$ then it remains a closed ideal in the unitization $\fu{A}$; therefore, by replacing $A$ with $\fu{A}$ if necessary, we may assume that $A$ has an approximate diagonal. Assume as before that $J$ has a right approximate identity bounded in norm by some constant $C>0$.

We now repeat the arguments above: approximate amenability ensures that we may choose $\Delta\in A\tp A$ satisfying \eqref{eq:approx-central} and \eqref{eq:pincer trick}, although we have no control on the norm of $\Delta$ itself. Nevertheless, writing $\Delta=\sum_{i=1}^m a_i\tp b_i$, we may choose an $f\in J$ with $\norm{f}\leq C$ such that $v\defeq \sum_{i=1}^m a_ifb_i$ satisfies $\norm{x\pi(\Delta)-xv} \leq\delta$ for all $x\in F$.
Hence, using \eqref{eq:three pieces}, we have $\norm{x-vx} \leq (C+2)\delta$ for all $x\in F$, which is enough to obtain a left approximate identity $(v_{F,\veps})$ for~$J$.

\end{section}




\vfill

\newcommand{\address}[1]{{\small#1.}}
\newcommand{\email}[1]{\texttt{#1}}

\noindent
\address{Yemon Choi,
Department of Mathematics and Statistics,
Lancaster University,
Lancaster LA1 4YF, United Kingdom} 

\noindent
\email{y.choi1@lancaster.ac.uk}

\end{document}